\documentclass[a4paper,twoside]{article}

\usepackage{epsfig}
\usepackage{subcaption}
\usepackage{calc}
\usepackage{romannum}
\usepackage{amssymb,amstext,amsmath,amsthm}
\usepackage{multicol}
\usepackage{pslatex}
\usepackage{natbib}
\let\bibhang\relax
\usepackage{apalike}
\usepackage{SCITEPRESS} 
\usepackage{multirow}
\usepackage{wrapfig}
\usepackage{graphicx,caption,varwidth}
\usepackage{blindtext}
\usepackage{rotating} % To display tables in landscape
\usepackage{booktabs}
\usepackage[linesnumbered,ruled,vlined]{algorithm2e}
\SetKwInput{KwInput}{Input}                % Set the Input
\SetKwInput{KwOutput}{Output}              % Set the Output
\usepackage{hyperref} 
\newcommand{\RomanNumeralCaps}[1]
{\MakeUppercase{\romannumeral #1}}
\begin{document}

\title{A LP relaxation based matheuristic for multi-objective integer programming}

 \author{\authorname{Duleabom An\sup{1}, Sophie N. Parragh\sup{1}\orcidAuthor{0000-0002-7428-9770},Markus Sinnl\sup{1}\orcidAuthor{0000-0003-1439-8702} and Fabien Tricoire\sup{2}\orcidAuthor{0000-0002-3700-5134}}
 \affiliation{\sup{1}Institute of 
 Production and Logistics Management, Johannes Kepler University Linz,\\ Altenberger Straße 69, 4040 Linz, Austria}
 \affiliation{\sup{2}Institute for Transport and Logistics Management,
 Vienna University of Economics and Business, \\Welthandelsplatz 1,
 1020 Vienna, Austria}
 \email{\{duleabom.an, sophie.parragh,\ markus.sinnl\}@jku.at, fabien.tricoire@wu.ac.at}
 }


\abstract{Motivated by the success of matheuristics in the single-objective domain, we propose a very simple linear programming-based matheuristic for three-objective binary integer programming.
To tackle the problem, we obtain lower bound sets by means of the vector linear programming solver \textit{Bensolve}. Then, simple heuristic approaches, such as rounding and path relinking, are applied to this lower bound set to obtain high-quality approximations of the optimal set of trade-off solutions. The proposed algorithm is compared to a recently suggested algorithm which is, to the best of our knowledge, the only existing matheuristic method for three-objective integer programming.
Computational experiments show that our method produces a better approximation of the true \textit{Pareto front} using significantly less time than the benchmark method on standard benchmark instances for the three-objective knapsack problem.}
\onecolumn \maketitle \normalsize \setcounter{footnote}{0} \vfill

\section{\uppercase{Introduction}}
\label{sec:introduction}
Many real-world optimisation problems involve multiple conflicting objectives, concerning, e.g., costs, environmental impact or service level, 
and can be formulated as integer linear programs.
%For example,\citet{fernandez2015triobjective} study, e.g., locating a public semi-desirable facility in a plane. The considered objectives are the minimisation of the 
%The authors introduce the mathematical model that minimises the 
%total distance between the users and the facility, and the total repulsion while aiming to distribute the repulsion equally.
\citet{kolli1999multiple}, for example, deal with facility location problems for a franchise company. When a new franchise is launched, the franchisee and the franchisor have conflicting objectives. To find the optimal location for the new store, the model maximises the number of customers while minimising conflict among existing franchises.
\citet{sawik2016bi} investigate vehicle routing problems occurring in the logistics of the food industry. The introduced model minimises both the total travel distance %of a vehicle 
and $CO_2$ emissions. %, considering constraints, e.g. start and endpoint of the route should be the same. 
A home care routing and scheduling problem is investigated by \citet{braekers2016bi}. Given a group of nurses and a number of care tasks at the patients' home locations, the model assigns tasks to each nurse while minimising operating cost and client inconvenience with respect to the timing of the visit and the nurse performing the task. 
\citet{kovacs2015multi} study a multi-objective consistent vehicle routing problem.
To provide consistent service in the routing industry, they maximise driver consistency and arrival time consistency while minimising total routing costs. % under restrictions such as time windows and the vehicle capacity.

The main goal in multi-objective (MO) optimisation  is to generate
the set of optimal trade-off solutions, known as \textit{efficient} (or \emph{Pareto optimal}) solutions.
Our study focuses on binary integer programming (IP) problems with three objectives. 

Motivated by the success of matheuristics in the single objective domain, we present a very simple matheuristic for three-objective integer programming.
To the best of our knowledge, only one other matheuristic for three-objective integer programming has been developed so far \citep{pal2019fpbh}. 
The algorithm we propose relies on a lower bound set which is obtained from the linear programming (LP) relaxation, using the vector linear programming solver $Bensolve$ \citep{lohne2017vector}.
The lower bound set is defined by its extreme points (and edges). Starting from those solutions which give rise to the extreme points, we apply \textit{rounding} in combination with \textit{path relinking}
to obtain high-quality approximations of the true \textit{Pareto frontier} ($PF$).\\
The contribution of this paper is twofold:
\begin{itemize}
\item[-] We show that, in the case of the three-objective assignment problem, since the extreme solutions which \textit{Bensolve} produces are integer, it already provides high-quality approximations of the true Pareto frontier, even without additional ingredients.
  \item[-] We propose the first LP relaxation-based matheuristic algorithm for three-objective integer programming, combining \textit{Bensolve} with \textit{rounding} and \emph{path relinking (PR)}
\end{itemize}

The remainder of the paper is structured as follows. Section \ref{sec::related work} briefly reviews existing work in multi-objective integer programming (MOIP) matheuristics. Section \ref{sec::preliminaries} provides
basic concepts and background for MOIP. 
The proposed algorithms are described in Section \ref{sec::algo frame} and empirically evaluated on the multi-objective assignment problem (MOAP) and multi-objective knapsack problem (MOKP) in Section \ref{sec::computational}. Finally, the paper concludes with a summary and suggestions for future work in Section \ref{sec::conclusion}.

\section{\uppercase{Related Work}} \label{sec::related work}
Over the past years a number of generic exact methods for solving MOIP have been proposed (see, e.g., \citet{mavrotas2013improved}, \citet{zhang2014simple}, \citet{kirlik2014new}, \citet{boland2017quadrant}). In spite of their popularity in the single objective domain, only comparably few contributions on matheuristc methods exist. 
%An large amount of research has been investigated to tackle MOIP. While the majority of studies have proposed exact approaches there are fewer contributions on matheuristics.

A heuristic and a matheuristic approach for binary
integer linear programming are proposed by \citet{soylu2015heuristic}. The methods are variants of variable neighbourhood search and local branching. Both algorithms collect segments of the $PF$ during the search
and combine them at the end. However, the proposed algorithms are not compared with a benchmark
method and struggle with generating a high-quality
approximation of the true $PF$ in quite a few problem
classes. To deal with bi-objective binary IP, \citet{leitner2016ilp} suggest an exact method and a %two-phase 
matheuristic framework.
In each phase, their heuristic obtains feasible solutions by fixing a large number of variables and reducing the associated feasible region, respectively. 
Based on the \textit{feasibility pump (FP)} idea introduced by \citet{fischetti2005feasibility},
\citet{pal2019feasibility} suggested a FP based heuristic for bi-objective IP and extended the method to higher dimensions in \citep{pal2019fpbh}. 
The authors suggest a two-stage heuristic algorithm
for MO mixed integer linear programming. The authors focus on finding diverse non-dominated points in the first stage by using the weighted sum and local search method.
In the second stage, they rely on local branching.
To the best of our knowledge, this is the only matheuristic method that deals with MOIP with more than two objectives. Therefore, we use it as a benchmark and refer to it as FPBH.

\section{\uppercase{Basic Concepts and Background}}\label{sec::preliminaries}
The problem we consider is a MOIP, with binary integer decision variables and three objectives. 
In the following, we state the MOIP in its general form. We present a unified view using minimisation objectives (any maximisation objective function can be
converted into a minimisation one by multiplying it
by -1).

\begin{equation}
\label{eq:MO}
    \begin{aligned}
        y = min\{C(x): Ax \geq b, x \in \{0,1\}^n \}, 
        \end{aligned}  \tag{MOIP}
\end{equation}\\
\noindent where $x_j$, $j=1,2,\dots,n$, is the vector of decision variables and $X:= \{Ax \geq b, x \in \{0,1\}^n\}$ is the feasible set. $C$ is a $p \times n$ objective function matrix where $c^k$, ($k=1,2,\dots,p$), is the $k^{th}$ row of $C$. In our case, $p=3$. Further, $y$ is a set of points in the objective space (criterion space), each of which corresponds to at least one solution vector $x \in X$. $A$ is an $m \times n$ constraint matrix and $b$ is the vector of right-hand-side values for these constraints.

\subsection{Pareto dominance and supported solutions}
In MO optimisation, the quality of a solution is determined by $Pareto$ dominance.
Suppose there are two solutions $x$ and $x^\prime$ of Problem (\ref{eq:MO}).
Then, $x$ \textit{dominates} $x^\prime$
if and only if $c^k(x) \leq c^k(x')$ for all $k \in \{1,\dots,p\}$ and $c^k(x) < c^k(x')$ for at least one $k$. If there does not exist any $x^\prime$ that dominates $x$, then $x$ is \textit{Pareto optimal}.
If $x^*$ is a \textit{Pareto optimal} (efficient) solution, then $y^*$ is a non-dominated point. The set of all non-dominated points is called the \textit{Pareto front}.\\
If a solution $x^*$ can be found by a convex combination of all objective functions, it is called a \textit{supported} efficient solution. Ohterwise, it is a \textit{non-supported efficient solution}.

\subsection{LP relaxation and bound set}
The notion of bound set was introduced by \citet{ehrgott2007bound}.
In the context of Problem (\ref{eq:MO}), a lower bound (LB) set and an upper bound (UB) set
%are denoted by $\underline{y}=(\underline{y}_1,\underline{y}_2,\underline{y}_3) $, and $\overline{y} = (\overline{y}_1,\overline{y}_2,\overline{y}_3) $, 
%such that $\underline{y} \leq y \leq \overline{y}$.
%Bounds sets 
provide information about the variable range that efficient solutions of the MOIP can attain. 
One common way to obtain a LB set for a minimisation problem is to solve the 
%solve a succession of weighted sum objectives of the  
LP relaxation of the original %for that
problem.
For constructing the LP relaxation of Problem (\ref{eq:MO}), the integrality conditions are removed, i.e. $0 \leq x_j \leq 1$, $j=1,\dots,n$ and $Bensolve$ can then be used to compute the LB set. In our heuristics, we will use the solutions associated with the LB set. In a slight abuse of notation, we will refer to these solutions as LB set.
% Given that the solutions of the LP relaxation problem are non-dominated by that of the IP, LP relaxation solutions can be a(n) LB/UB set of (minimisation/maximisation) IP,
%In our context, we obtain a LB set for minimisation problems. 

\subsection{Benchmark problems}\label{sec::benchmark prob}
The multi-objective assignment problem (MOAP) and the multi-objective knapsack problem (MOKP) are commonly used benchmark problems in MO. We use standard benchmark instances of these two problems for our computational experiments.

\subsubsection{Multi-objective assignment problem}
The well-known assignment problem is a type of transportation problem: a certain number of tasks $l \in \{1,\dots,n\}$ and agents $r \in \{1,\dots,n\}$ are given. The decision variable $x_{rl}$ denotes whether task $l$ is assigned to the agent $r$ ($x_{rl}$=1) or not ($x_{rl}$=0). When a task is allocated to an agent, the corresponding non-negative costs $c^1_{rl},\dots,c^p_{rl}$ are incurred. The MOAP can be stated as follows:\\

\begin{equation}
\label{eq:AP1}
\begin{aligned}
    &\min \sum_{r=1}^{n}\sum_{l=1}^{n}{c_{rl}^{j}x_{rl}} & j=1,\ldots,p \\
    % && r,l= 1,\dots,n\\
\end{aligned}  
\end{equation}
\begin{equation}
\label{eq:AP(2)}
\begin{aligned}
    &\textrm{s.t.} \sum_{l=1}^{n}{x_{rl}} = 1 & r = 1,\dots,n \\
\end{aligned}  
\end{equation}
\begin{equation}
\label{eq:AP(3)}
\begin{aligned}
    &\sum_{r=1}^{n}{x_{rl}} = 1 & l = 1,\dots,n \\
    \end{aligned}   
\end{equation}\\
\begin{equation}
\label{eq:AP4}
\begin{aligned}
    &x_{rl}\in \{0,1\}  & r,l=1,\ldots,n.
\end{aligned}  
\end{equation}\\
The objective of the MOAP is to find an optimal assignment of all tasks to agents while minimising $p$ cost functions (\ref{eq:AP1}). Equation (\ref{eq:AP(2)}) ensures that each agent is assigned to only one task. Equation (\ref{eq:AP(3)}) limits each task to be assigned to one agent only.

It is well-known that the constraint matrix of AP is totally unimodular. According to the property, we can conclude that every solution of an LP relaxation is naturally an integer vector.
 %However, there may be an optimal solution  of LP relaxation that is not a vertex of the feasible region.
% Although $Bensolve$ only computes supported efficient solutions, these solutions already 
Although it is restricted to supported efficient solutions as $Bensolve$ is based on Benson's outer approximation, it already provides a high-quality approximation of $PF$. The corresponding computational results are given in Table \ref{tb:AP} in Section \ref{sec::computational}.
In conclusion, the MOAP may not be a suitable benchmark for MOIP heuristics.
%, although it has been used as a test instance in studies of MOIP, e.g. \citep{pal2019fpbh}. 

 \subsubsection{Multi-objective knapsack problem}
In the multi-objective knapsack problem, a set of items is given, each with a certain weight $w_r$, and we must select a subset of these items such that the total weight does not exceed a given capacity $W$. The decision variable $x_r$ denotes whether item $r$ is selected for the knapsack ($x_r$=1) or not ($x_r$=0). Each item $r$ also has profits $v_r^1 \dots v_r^p$. Here, $W$, $v_r$, and $w_r$ are non-negative integer values. 
The MOKP model is stated as follows:\\

\begin{equation}
\label{eq:KP(1)}
\begin{aligned}
    &&\max \sum_{r=1}^{n}{v_{r}^{j}x_{r}} && j=1,\ldots,p \\
\end{aligned}  
\end{equation}
\begin{equation}
\label{eq:KP(2)}
\begin{aligned}
    &\textrm{s.t.} \sum_{r=1}^{n}{w_{r}x_{r}} \leq W \\
\end{aligned}  
\end{equation}
\begin{equation}
\label{eq:KP(3)}
\begin{aligned}
    &x_{r}\in \{0,1\}  & r=1,\ldots,n 
\end{aligned}  
\end{equation}\\
Equation (\ref{eq:KP(1)}) denotes the objective functions maximising the $p$ total profits of selected items. Equation (\ref{eq:KP(2)}) is the capacity constraint. The total weight of the items placed in the knapsack cannot exceed the given capacity.\\ 
In this paper, we convert a maximisation objective function into a minimisation one by multiplying it by~$-1$.
%the objective function by 
 % for unification. %unified to% minimisation problems, consequently the bound set provided by Bensolve is always an LB set. 

\section{\uppercase{LP relaxation-based matheuristic}}
\label{sec::algo frame}
This section provides the overview of the proposed algorithm and its main ingredients.
%explanations of simple heuristic approaches applied to LB sets. 

\subsection{Algorithmic framework\label{algo:ben+pr}}
The proposed matheuristic algorithm follows a two-stage approach. 
At the first stage, we obtain LB sets and adapt them to feasible integer rounded sets.
A variant of the weighted sum method by \citet{ozpeynirci2010exact}, \textit{Bensolve} by \citet{lohne2017vector} and \textit{Inner approximation} by \citet{csirmaz2020inner} have been investigated to find LB sets first. Among the three methods, $Bensolve$ shows competitive performance in the experiment and produces all the required bound set information we need. The performance comparison can be found in the Appendix.
In the context of MOKP, once LB sets ($L$) are obtained, we round down the fractional variable to produce a feasible solution. These are referred to as integer rounded (IR) sets. At the second stage, $PR$ is employed to improve the solution quality. 
Until reaching the iteration limit, the $PR$ process repeats. If $PR$ finds a new solution, it is stored in the archive $candX$.
%As computational time and the number of solutions grow proportionally to the number of iterations, a decision-maker can adjust it depending on a condition such as a time constraint. 
Based on initial experiments, we set the limit of iterations of $PR$ to \textit{the number of solutions of the IR set times 50}.
Since the dominance relation is not checked in $\Tilde{X}$ during the search, dominated solutions are filtered after the algorithm terminates, and the set of integer feasible solutions $X$ is returned.
Algorithm \ref{algo:framework} describes the general framework of the proposed algorithm.
\begin{algorithm}
 \KwInput{$L$: points describing an LB set}
 $candX$: an archive of newly found feasible IP solutions by $PR$\\
 $\Tilde{X}$: an empty list\\
%   \textit{ratio\_table}: $m \times n$ matrix, where $m$ and $n$ are the number elements in $ND$ and objectives ($p$=3), respectively.\\
  \vskip 2mm
  $IR \leftarrow$ \textit{RoundingDown}($L$); \\
%   $candX \leftarrow \emptyset$;\\

  $i=0$;\\
  \For{i$<$ iteration limit}
  {
   $candX \leftarrow$\textit{PathRelinking}($IR,S_i,S_g,IGPair$);\\
   Update $candX$\\
    % $\Tilde{X} \leftarrow \Tilde{X} \cup candX$;\\
    % $IGPair \leftarrow [S_i,S_g]$;\\
    $i=i+1$;
  }
  $\Tilde{X} \leftarrow$ \textit{DominanceCheck}($candX$);
  \vskip 2mm
 \KwOutput{$\Tilde{X}$}
 \caption{The LP relaxation-based matheuristic framework}
 \label{algo:framework}
\end{algorithm}

% \subsection{Rounding}
% Solutions to the LP relaxation of the KP always feature one fractional variable, the one item which does not fit entirely into the knapsack. Rounding this variable value "down" to zero will result in a feasible integer solution. We use this property to obtain feasible integer solutions from the LB set produced by \emph{Bensolve}. We call this set integer rounded (IR) set.
%In the context of KP, items with a higher \textit{profit/weight} ratio are added into the knapsack first to maximise the total profit. Accordingly, in a relaxation problem of KP, the decision variable of the lastly added item usually has a fractional value. By rounding down this fractional value, the solution set automatically becomes feasible. This is how the IR set is produced.
% Based on this property, we round down the fractional decision variable in the LP solution that represents the lastly added item into the knapsack. As a consequence, the adapted solution becomes an integer feasible solution and is referred to as integer rounded (IR) sets. This is the first stage of the proposed algorithm obtaining IR sets.

\subsection{Path relinking}
Although IR sets obtained from \textit{rounding} are feasible, they are not necessarily of high-quality. 
Therefore, we employ $PR$ to increase the number of solutions and improve the quality of the approximate \textit{PF}.
% The set of rounded down LP solutions (integer feasible solutions), denoted $IR$, is a reference set. 
% Before reaching the iteration limit, the algorithm
% repeats the whole $PR$ process. 
% In the second stage, we extend our heuristic  The employed approach is $PR$. 
$PR$ was originally introduced by \citet{glover1997tabu}. 
The main idea of it is that there should be common characteristics among high-quality solutions. The method 
%exploits the common attributes for producing 
produces new solutions by exploring solution "paths" between pairs of known solutions.
To generate a new solution, $PR$ chooses two solutions from a set of initial solutions; an initiating solution ($S_i$) and a guiding solution ($S_g$) represent the starting and ending points of the path, respectively. Then it explores the trajectories that link the pair of solutions in a neighbourhood space.  
In each step, a certain number of intermediate paths, called a neighbourhood, are created. 
For the next step, the new initiating solution is selected in there.
% From the neighbourhood, one solution is selected and set to the new initiating solution for the next step. 
As the search proceeds, more attributes of the guiding solutions are gradually passed on to intermediate solutions.
The search continues until the initiating solution becomes identical to the guiding solution.

The method has been applied successfully to MO problems such as the travelling salesman problem \citep{jaszkiewicz2005path}, the dial-a-ride problem \citep{parragh2009heuristic}, and a school bus routing problem \citep{de2017multi}. It also produces high-quality solutions on large MOKP instances within a reasonable timeframe \citep{beausoleil2008multi,marti2015multiobjective}. 
For these reasons, we combine $PR$ with the first stage approach.

For illustration purposes, we provide a short example describing the $PR$ process on a MOKP with three objectives. Suppose we have a four-item MOKP in which the initiating solution is (0 0 1 0) and the guiding solution is (1 1 0 0). In this case, the only selected item in common is the fourth item. Therefore, three intermediate paths are created by the following rules. To create one path, one variable value in the initiating solution is switched. To be specific, if the item is placed in the knapsack and its value is 1, it is taken out of the knapsack and the value changes to 0. If the item is not chosen and the value equals 0, it is placed in the knapsack and the value changes to 1. This procedure applies to all the different variable values. The outcome of the process is seen in the first row (Neighbourhood) of Table \ref{table:PR}.
When the $P$ matrix (\ref{eq::profits}) represents the profits of the three objectives, the sets of total profits that correspond to the first neighbourhood are [7,6,8], [5,4,6], and [0,0,0]. Since the dominance relation among the solution points is clear, the first path (1 0 1 0) is selected as the new initiating solution (marked with * in Table \ref{table:PR}). This process repeats until the initiating solution becomes identical to the guiding solution. Table \ref{table:PR} illustrates the transforming process of $PR$ in this example.
\begin{equation}\label{eq::profits}
P = 
\begin{pmatrix}
4 & 2 & 3 & 6\\
5 & 3 & 1 & 8\\
6 & 4 & 2 & 7
\end{pmatrix}
\end{equation}\\

\begin{table}[ht]
\begin{tabular}{|c|c|c|c|c|}
  \hline
  Initiating & Guiding & \multicolumn{3}{|c|}{Neighbourhood} \\
  \hline
  0 0 1 \textbf{0} & 1 1 0 \textbf{0} &1 0 1 0* & 0 1 1 0 & 0 0 0 0  \\
  \cline{1-5}
  \textbf{1} 0 1 \textbf{0}&\textbf{1} 1 0 \textbf{0} &  1 1 1 0* & 1 0 0 0  &\multicolumn{1}{c}{}\\ 
  \cline{1-4}
  \textbf{1} \textbf{1} 1 \textbf{0}&\textbf{1} \textbf{1} 0 0 & 1 1 0 0* &\multicolumn{2}{c}{}\\
  \cline{1-3} 
  \textbf{1 1 0 0} & \textbf{1 1 0 0} &  \multicolumn{3}{c}{}\\
  \cline{1-2}
\end{tabular}
\caption{\textit{Path relinking} procedure}\label{table:PR} \centering
\end{table}
In our $PR$ process, an infeasible intermediate solution is not stored in the solution archive, but still can be a new initiating solution. This is for letting the algorithm explore a broader search region so that it possibly finds diverse solutions.

Depending on the specific problem, the following components of $PR$ can be designed differently in the algorithm.
\begin{itemize}
\setlength\itemsep{0em}
    \item[1)] Building an initial solution set
    \item[2)] Selecting $S_i$ and $S_g$
    \item[3)] Generating intermediate paths (a neighbourhood)
    \item[4)] Choosing the next initiating solution
\end{itemize}
In this study:
\begin{itemize}
\setlength\itemsep{0em}
    \item[1)] We use IR sets as the initial solution set.
    \item[2)] $S_i$ and $S_g$ are chosen 
    either randomly or based on $similarity$ between them.
    Since the performance of $PR$ can vary depending on the choice of the $S_i$-$S_g$ pair, we investigate different ways, and call \textit{SelectionRule\RomanNumeralCaps{1}}.
    The first approach is to select two random solutions from $IR$. For the second and third methods, $S_i$ is chosen at random in $IR$ first, then the $similarity$ between $S_i$ and all the other solutions in $IR$ is calculated. The $similarity$ is measured by counting the number of variables which have equal values. Thus, higher values imply a greater $similarity$. 
    The second method selects $S_g$ by finding the most similar solution to $S_i$. The third approach finds the most different solution to $S_i$ for $S_g$. When $S_i$ and $S_g$ are similar, fewer intermediate paths are likely to be created. However, when they are different, the $PR$ process can take longer as it explores relatively diverse intermediate paths.  

    \item[3)] A new neighbourhood is generated as many times as the number of different variable values.
    \item[4)] The next $S_i$ is determined either by dominance relation based analysis or randomly.
    % In the first case (acceptance), we analyse the dominance relation among intermediate solutions and select the best solution as $S_i$. For this, \textit{SelectionRule\RomanNumeralCaps{2}} is introduced. If one non-dominated solution exists in the neighbourhood, that becomes the next initial solution. If there are mutually non-dominating solutions, we check the \textit{improved ratio} of each solution point. The approach for finding out the most-improved point is explained in detail in Section \ref{sec::improved}. For the second case (rejection), the analysis is not needed, instead one random solution from the neighbourhood is chosen.
\end{itemize}
The description of the entire $PR$ implementation for the proposed algorithm is given in Algorithm \ref{algo:BenPR}.

The algorithm maintains two archives $candX$ and $IGPair$.
% which stores newly found solutions.
% two archives 
Each archive stores newly found solutions and used $S_i$-$S_g$ pairs during the $PR$ process.
Under the iteration limit, the whole process repeats.
% Since the performance of $PR$ can vary depending on the choice of the $S_i$-$S_g$ pair, we investigate different ways, and call \textit{SelectionRule\RomanNumeralCaps{1}}. The first approach is to select two random solutions from $IR$. For the second and third methods, $S_i$ is chosen at random in $IR$ first, then the $similarity$ between $S_i$ and all the other solutions in $IR$ is calculated. The $similarity$ is measured by counting the number of variables which have equal values. Thus, higher values imply a greater $similarity$. The second method selects $S_g$ by finding the most similar solution to $S_i$. The third approach finds the most different solution to $S_i$ for $S_g$. When $S_i$ and $S_g$ are similar, fewer intermediate paths are likely to be created. However, when they are different, the $PR$ process can take longer as it explores relatively diverse intermediate paths.\\
% Once $S_i$ and $S_g$ are chosen by  \textit{SelectionRule\RomanNumeralCaps{1}} (line 4),
% $PR$ continues until one of the following two termination conditions is met (line 5):
% The description of the entire $PR$ implementation for the proposed algorithm is given in Algorithm \ref{algo:BenPR}.
% The algorithm maintains two archives $candX$ and $IGPair$.
% Each archive stores newly found solutions and used $S_i$-$S_g$ pairs during the $PR$ process. Under the iteration limit, the whole process repeats.
Firstly, $S_i$ and $S_g$ are chosen by  \textit{SelectionRule\RomanNumeralCaps{1}} (line 7).
A $PR$ iteration ends when one of the following two termination conditions is met (line 8);
\begin{itemize}
\item[-] The initial and guiding solutions are identical.
\item[-] The pair of $S_i$ and $S_g$ has already been chosen.
\end{itemize}
To decide the number of intermediate paths, we need to find the positions of the differing variable values in $S_i$ and $S_g$. 
They are stored in $\Delta$items as indices (lines 9-10).
Once $\Delta$items is defined, we generate the neighbourhood by the following rules: For each index in $\Delta$items, the variable values in $S_i$ change one by one (line 11). Let us suppose items 1 and 5 have a different value, $\Delta$items = $\{1,5\}$. If item 1 is selected in $S_i$  (variable value=1), then it is taken out of it. If it is out of the knapsack (variable value=0), then it is added to it. The same process is conducted on item 5. 
Once intermediate paths are built, we select the next $S_i$.
% analyse them to select the best path (line 12). 
In order to avoid local optima, the best neighbour is not systematically selected. Rather, it is selected with a certain probability, set to 0.7 after initial experiments. Otherwise, another neighbour is selected randomly.
In the case where we do want to select the best neighbour, we analyse the dominance relation among intermediate solutions and select the best solution as $S_i$ (line 14). For this, \textit{SelectionRule\RomanNumeralCaps{2}} is introduced. If one non-dominated solution exists in the neighbourhood, it becomes the next initiating solution. If there are mutually non-dominating solutions, we check the \textit{improved ratio} of each solution point. The approach for finding the most-improved point is explained in detail in Section~\ref{sec::improved}. For the case where we do not want to select the best neighbour, the analysis is not needed and simply one random solution is selected from the neighbourhood (lines 15-16).
% In the acceptance case, if there is a dominating solution in the neighbourhood, that becomes the new initiating solution (lines 14-15). Otherwise, where there are mutually non-dominating solutions, $S_i$ is selected by \textit{SelectionRule\RomanNumeralCaps{2}} (lines 16-17). 
% When the chosen best solution is rejected, one random path from the neighbourhood is selected and set to the new $S_i$ (lines 18-19). 
Before moving on to the next iteration, the algorithm checks the feasibility of $S_i$ (line 18). 
If $S_i$ is feasible and not included in $IR$, it is stored in the archive $candX$ (line 19). The infeasible $S_i$ is not archived, though, it is still used for the next iterations as it might help to find additional feasible solutions in yet unexplored parts of the search region. A newly found feasible solution is added to $IR$ (line 20).
To prevent the case in which the same pair of initial and guiding solution is chosen, the current $[S_i,S_g]$ pair is archived after every $PR$ iteration (line 21). 
% The process repeats until $S_i$ becomes equal to $S_g$.
% Since we do not check the dominance relation in the archive $candX$ during the search, 
After post-processing that filters dominated solutions (line 23),
The algorithm returns the set of integer feasible solutions $\Tilde{X}$.
% is conducted after the algorithm terminates.   

\begin{algorithm}
 \KwInput{$IR$}
  $candX$: an archive of newly found feasible IP solutions\\
  $IGPair$: an archive of $S_i$-$S_g$ pairs\\
  $\Tilde{X}$: an empty list\\
  \vskip 1mm
  $candX \leftarrow \emptyset$;   $IGPair\leftarrow \emptyset$; \\
  $i$=0;\\
  \For{i $<$ iteration limit}
  {
  Select $S_i$,$S_g$ from $IR$ following \textit{SelectionRule\RomanNumeralCaps{1}};\\
    \vskip 1mm
     \While{$S_i \neq S_g$ and [$S_i,S_g$] $\notin$ IGPair}{
      $\Delta$ items $\leftarrow$  index set of different variable values;\\
      $n \leftarrow$ \#indices in $\Delta items$;\\
     Create $n$ neighbourhood;  \\
    %  \For{i $\in \Delta$ items}{
    %     % Switch the state of $i$ item in $S_i$ and archive in \textit{neighbourhood}\\
    %     \eIf {$S_i$[$\Delta$ items[$i$]] == 0}
    %     {$S_i$[$\Delta$ items[$i$]] = 1; }
    %     {$S_i$[$\Delta$ items[$i$]] == 0;}
    %  }
      The best move analysis: \\
      \eIf{rand() $<$ 0.7}{
            $S_i$ is chosen by \textit{SelectionRule\RomanNumeralCaps{2}}
            % Check dominance relation among $neighbourhood$;\\
            % \eIf{one strongly dominating solution $ new S_i$ exists}
            % {$S_i \leftarrow newS_i$;}
            % { }; 
            % either
            %  Select one of the dominating solutions $new S_i$; \\ $S_i \leftarrow newS_i$;\\
            % or
            %  $S_i \leftarrow ImprovedND$;
             
      }
      {
      $S_i \leftarrow$ randomly choose one solution from the neighbourhood; 
      }
      Feasibility check of $S_i$;\\
      \If{$S_i$ \,is\, feasible and $S_i \notin IR$}{
      $candX \leftarrow S_i$;\\
    %   $IR \leftarrow IR \cup S_i$;\\
      Update $IR$ by adding $S_i$;  
      }
     }
    $IGPair \leftarrow [S_i,S_g]$;\\
    $i=i+1$;\\
 }
  $\Tilde{X} \leftarrow$ \textit{DominanceCheck}($candX$);\\
  \KwOutput{$\Tilde{X}$}
 \caption{LP relaxation-based matheuristic}
 \label{algo:BenPR}
\end{algorithm}

\subsubsection{$ImprovedND$ operation}\label{sec::improved}
The \textit{ImprovedND} operation figures out which path shows the biggest improvement compared to the current solution $S_i$. Algorithm \ref{algo:improved} shows a precise description of the operation.

Let $ND$ be a set of non-dominated intermediate points and $objS_i$ be the objective values of $S_i$. To record the improvement of each non-dominated point, we create two matrices and one list. The \textit{ratio\_table} stores the ratio of a non-dominated point to the current point $objS_i$, for each objective (lines 4-6). Afterwards, we assign a rank to each column of the \textit{ratio\_table} and enter rankings into the \textit{rank\_table} (lines 7-10). The rankings of each non-dominated point are added up (lines 11-12) and stored in \textit{ND\_degree}. The non-dominated path with the highest degree is set to $S_i$ (line 13).
\begin{algorithm}
 \KwInput{$objS_i$, $ND$}
  \textit{ratio\_table}: $|ND| \times p$ matrix\\
  \textit{rank\_table}: $|ND| \times p$ matrix\\
  \textit{ND\_degree} : list of size $|ND|$
  \vskip 2mm
 
  \For{i=1,...,$|ND|$ }
  {
    \For{j=1,...,$p$}
        {
        \textit{ratio\_table}[$i$][$j$] = $\cfrac{ND[i][j]}{objS_i[j]}$ 
        }
    }
  \vskip 2mm
  
    \For{i=1,...,$|ND|$}
    {
        \For{j=1,...,$p$}
            {
            \textit{rank\_table}[$i$][$j$] = rank of $ND$[$i$][$j$] \\ in $j^{th}$ column
            }
    }
    \vskip 2mm
    \For{i=1,...,$|ND|$}
    {
    \textit{ND\_degree}[$i$] $\leftarrow$ sum(\textit{rank\_table} $i^{th}$ column)\\
    }
    $S_i \leftarrow$ $ND$ with the highest degree 

  \KwOutput{$S_i$}
 \caption{$ImprovedND$}
 \label{algo:improved}
\end{algorithm}

\section{\uppercase{Computational Experiments}}\label{sec::computational}
% We evaluate the performance of two versions of our heuristic method, especially focusing on the comparison with the benchmark, FPBH \citet{pal2019fpbh}. 
We evaluate the performance of the different versions of our heuristic method, focusing on the comparison with FPBH \citep{pal2019fpbh}. The variant with \textit{rounding down} is named $RD$. The following $PR$ variants are tested. 
The first three $PR$ variants randomly choose the next initiating solution during the iterations.
\begin{itemize}
 \setlength\itemsep{-0.1em}
\item \textit{PRrand}: This version randomly selects both the initiating and guiding solutions from the initial solution set.
\item $PR$sim: This version randomly chooses the initiating solution, then finds the most similar solution among the remaining solution set for the guiding solution.
\item $PR$dif: This version finds the most different solution to the already chosen initiating solution for the guiding solution.
\end{itemize}
The three other variants consider the improvement among mutually non-dominating intermediate paths to choose the next initiating solution.
\begin{itemize}
\setlength\itemsep{-0.1em}
\item \textit{PI}: This version selects the most improved intermediate solution as the next initiating solution.
\item $PI$sim: This variant uses \textit{ImprovedND} within $PR$sim.
\item $PI$dif: This variant uses \textit{ImprovedND} in $PR$dif.
\end{itemize}

The proposed algorithms use \textit{Bensolve} by \citet{lohne2017vector} to obtain the bound sets. The heuristic integration (\textit{rounding down} and $PR$), is implemented in Julia. For the benchmark algorithm, we used the Julia implementation of FPBH (with the default setting) which is publicly available at \url{https://github.com/aritrasep/FPBH.jl}. All experiments of matheuristics are carried out on Intel® Core™ i5-8250U CPU running at 1.60GHz with 16GB RAM. The exact MOIP solver proposed by \citet{kirlik2014new} (KS) is also used in the experiment to obtain the true \textit{PF} for comparison purpose. KS is run on Quad-core X5570 Xeon CPUs @2.93GHz with 48GB RAM. The KS results are for reference only (i.e. not for benchmarking).

The test instances we use are the same ones on which FPBH is tested. The instances were generated by \citet{kirlik2014new} and are publicly available at \url{http://home.ku.edu.tr/~moolibrary/}. 
Each problem class has 100 instances divided into 10 subclasses, each of which contains 10 instances. MOAP instances are formed in the number of tasks (to be assigned) which varies from 5 to 50 in increments of 5. MOKP instances are classified by the number of items which varies from 10 to 100 in increments of 10. 

\subsection{Performance measure: Hypervolume indicator}
One widely used indicator to measure the quality of a solution set in MO optimisation is the hypervolume (HV) indicator. HV measures the volume of the dominated space of all the solutions contained in a solution set. To calculate the dominated space, a reference point must be used. Usually, a reference point is the “worst possible” point in the objective space. 
In this study, all the HV values are calculated with normalised objective values. \\
Let $Y^{k}_N$ be a set of $k^{th}$ objective values of the true $PF$ and $y\in \mathbb{R}^{p}$ be an arbitrary non-dominated point obtained from a heuristic algorithm. Then, the normalised values of the obtained point are:
\begin{equation*}
    \begin{aligned}
        \frac{y^k-min(Y^{k}_N) }{max(Y^{k}_N)-min(Y^{k}_N)} && k=1,\dots,p.
    \end{aligned} 
\end{equation*}
As all non-dominated points are normalised, their values exist in [0,1]. Therefore, the reference point is (1,1,1). Higher HV values indicate a better approximation. We used the publicly available HV computing program provided by \citet{fonseca2006improved} at \url{http://lopez-ibanez.eu/hypervolume#intro} to obtain HV values.

\subsection{Results and discussion}\label{results}
We report the following results of each algorithm: the number for solutions ($|Y|$), CPU time (sec), HV value, and HV as a percentage of the HV indicator value for the exact non-dominated set as provided by \citet{kirlik2014new}. 
% that demonstrates how far the algorithm reaches the exact \textit{PF}. 
% There is no time limit, as all the algorithms finish searching once they meet the termination condition.
All the figures are average results over 10 test instances. The figures of $PR$ variants and FPBH are averaged over 10 runs for each instance because they have random components. The results of the experiments on MOAP and MOKP are reported in Tables \ref{tb:AP}-\ref{tb:KPhv}.

$Bensolve$ already finds integer feasible solutions for all the MOAP instances. In addition, it shows better performance than FPBH in all the subclasses. Therefore, we do not apply our matheuristic to MOAP instances.

% \begin{sidewaystable}
% \centering
\begin{table*}[!ht]
\makebox[\linewidth]{
 \begin{tabular}{crrrrrrrrrrrrrrr} \toprule %\hline
    n  &\multicolumn{3}{c}{$|Y|$} & &\multicolumn{3}{c}{CPUtime(sec)} & & \multicolumn{3}{c}{HV} & &\multicolumn{2}{c}{HV(\%)} \\ \cline{2-4} \cline{6-8} \cline{10-12} \cline{14-15}
    &  KS* & FPBH & $Bensolve$  & & KS& FPBH & $Bensolve$ & &KS & FPBH & $Bensolve$& & FPBH & $Bensolve$ \\ \midrule%\hline
    5  & 14.1 & 6.7 & \textbf{7.5} && 0.11 & 0.05 & \textbf{0.001} && 6.77 & 6.60 & \textbf{6.71} && 97.49 &\textbf{99.11}\\ 
    10  & 176.8 & 21.0 & \textbf{39.0} && 10.71 & 0.28 & \textbf{0.019} && 7.23 & 6.92 & \textbf{7.18} && 95.71 &\textbf{99.31}\\  
    15  & 674.9 & 40.4 & \textbf{83.1} && 92.51 & 1.17 & \textbf{0.068} && 7.33 & 6.96 & \textbf{7.29} && 94.95 &\textbf{99.45}\\ 
    20  & 1860.5 & 62.5 & \textbf{161.3} && 359.07 & 2.04 & \textbf{0.222} && 7.40 & 6.99 & \textbf{7.37} && 94.46 &\textbf{99.59}\\ 
    25  & 3567.8 & 90.6 & \textbf{253.1} && 872.19 & 4.99 & \textbf{0.596} && 7.46 & 7.05 & \textbf{7.44} && 94.50 &\textbf{99.73}\\
    30  & 6181.3 & 140.0 & \textbf{379.4} && 1859.74 & 15.59 & \textbf{1.157} && 7.48 & 7.04 & \textbf{7.46} && 94.12 &\textbf{99.73}\\
    35  & 8972.3 & 163.2 & \textbf{501.4} && 3285.57 & 26.39 & \textbf{2.082} && 7.50 & 7.06 & \textbf{7.48} && 94.13 &\textbf{99.73}\\ 
    40  & 14679.7 & 242.9 & \textbf{699.1} && 6425.98 & 57.95 & \textbf{3.824} && 7.53 & 7.09 & \textbf{7.51} && 94.16 &\textbf{99.73}\\
    45  & 17702.2 & 238.6 & \textbf{838.0} && 9239.01 & 82.16 & \textbf{6.097} && 7.56 & 7.10 & \textbf{7.54} && 93.92 &\textbf{99.74}\\
    50  & 24916.8 & 337.5 & \textbf{1034.8} && 15814.82 & 119.52 & \textbf{9.739} && 7.58 & 7.12 & \textbf{7.56} && 93.93 &\textbf{99.74}\\
    \bottomrule
\end{tabular}
}
\caption{Comparing algorithm performance on MOAP for $p$=3, * indicates optimal values, best heuristic values are in bold.} 
\label{tb:AP}
\end{table*}
% \end{sidewaystable}

\begin{table*}[!ht]
% \centering
\makebox[\linewidth]{
 \begin{tabular}{crrrrrrrrrrr} \toprule %\hline
    &&&&&\multicolumn{6}{c}{$PR$ variants} & \\ \cline{6-11} 
    &n&  KS* & FPBH & $RD$ & $PRrand$ & $PR$sim& $PR$dif & $PI$ &$PI$sim &$PI$dif & \\ \midrule%\hline
    &   10 &9.8&	5.1&    4.3	&   5.7&	4.8&	4.8&	\textbf{6.3}&	4.6&	4.6& \\
    &   20 &38.0&	18.1&	10.7&	22.4&	20.7&	19.0&	\textbf{26.1}&	20.4&	18.6& \\
    &   30 &115.8&	43.2&	20.7&	47.9&	46.9&	45.6&	\textbf{58.6}&	46.3&	45.0& \\
    &	40&	311.2&	95.7&	33.8&	95.8&	96.4&	91.3&	\textbf{122.0}&	97.8&	93.2& \\
    &	50&	444.2&	111.8&	41.7&	119.1&	118.1&	114.9&	\textbf{143.5}&	118.4&	112.6&	\\
    &	60&	917.1&	195.1&	71.5&	209.1&	209.8&	203.6&	\textbf{266.2}&	207.8&	208.4&	\\
    &	70&	1643.4&	348.2&	90.2&	263.0&	264.4&	259.8&	\textbf{353.4}&	271.3&	264.9&	\\
    &	80&	2295.8&	\textbf{439.0}&	113.1&	305.1&	301.1&	310.1&	399.9&	313.7&	309.2&	\\
    &	90&	3107.8&	\textbf{501.9}&	130.6&	322.7&	319.4&	327.0&	412.7&	338.7&	343.3&	\\
    &	100&5849.0&	\textbf{919.2}&	176.7&	442.8&	453.0&	439.8&	581.3&	469.4&	471.9&	\\
    \bottomrule
\end{tabular}
}
\caption{Comparing $|Y|$ of algorithms on MOKP for $p$=3, * indicates optimal values, best heuristic values are in bold.} 
\label{tb:KPsolution}
\end{table*}

\begin{table*}[!ht]
% \centering
\makebox[\linewidth]{
 \begin{tabular}{crrrrrrrrrrr} \toprule %\hline
    &&&&&\multicolumn{6}{c}{$PR$ variants} & \\ \cline{6-11} 
    &n&  KS* & FPBH & $RD$ & $PRrand$ & $PR$sim& $PR$dif & $PI$ &$PI$sim &$PI$dif & \\ \midrule%\hline
    &	10	&	0.140	&	0.023	&	\textbf{0.001}	&	0.002	&	0.004	&	0.003	&	0.006	&	0.002	&	0.002	\\
    &	20	&	1.030	&	0.080	&	\textbf{0.004}	&	0.020	&	0.056	&	0.044	&	0.067	&	0.051	&	0.036	\\
    &	30	&	5.540	&	0.324	&	\textbf{0.008}	&	0.085	&	0.314	&	0.281	&	0.307	&	0.312	&	0.268	\\
    &	40	&	23.23	&	1.071	&	\textbf{0.013}	&	0.224	&	0.949	&	0.896	&	0.824	&	0.885	&	0.763	\\
    &	50	&	40.07	&	1.941	&	\textbf{0.019}	&	0.379	&	1.752	&	1.629	&	1.517	&	1.684	&	1.416	\\
    &	60	&	116.0	&	5.332	&	\textbf{0.041}	&	1.105	&	5.490	&	5.095	&	4.872	&	5.445	&	4.686	\\
    &	70	&	283.5	&	12.68	&	\textbf{0.054}	&	2.118	&	10.69	&	9.818	&	10.72	&	10.64	&	9.118	\\
    &	80	&	440.0	&	20.77	&	\textbf{0.079}	&	3.533	&	18.18	&	16.21	&	17.49	&	17.98	&	15.31	\\
    &	90	&	833.9 &	42.17	&	\textbf{0.102}	&	6.121	&	28.89	&	26.50	&	30.04	&	28.18	&	24.41	\\
    &	100	&	2478.4	&	82.54	&	\textbf{0.129}	&	11.23	&	59.78	&	57.50	&	66.97	&	60.57	&	53.24	\\
    \bottomrule
\end{tabular}
}
\caption{Comparing CPU time (sec) of algorithms on MOKP for $p$=3, * indicates optimal values, best heuristic values are in bold.} 
\label{tb:KPcpu}
\end{table*}

%%%%%%%%%%%%%%%%%%%%%%for presentation 
% \begin{table*}[!ht]
% % \centering
% \makebox[\linewidth]{
%  \begin{tabular}{crrrrrrrrrrr} \toprule %\hline
%     && \multicolumn{4}{c}{CPUtime (sec)} && \multicolumn{3}{c}{HV(\%)}\\ \cline{3-6} \cline{8-10} 
%   &	n	&	KS*	&	FPBH	&	$RD$	&	$PRvar$	&&	FPBH	&	$RD$	&	$PRvar$ \\ \midrule%\hline
% &	 10 	&	 0.140 	&	 0.023 	&	 \textbf{0.001} 	&	0.006	&	&	\textbf{95.44}	&	89.82	&	91.64\\
% &	     20 	&	 1.030 	&	 0.080 	&	 \textbf{0.004} 	&	0.067	&	&	96.13	&	94.84	&	\textbf{97.13}\\
% &	     30 	&	 5.540 	&	 0.324 	&	 \textbf{0.008} 	&	0.307	&	&	95.98	&	96.53	&	\textbf{97.92}\\
% &	     40 	&	 23.23 	&	 1.071 	&	 \textbf{0.013} 	&	0.824	&	&	96.08	&	97.34	&	\textbf{98.46}\\
% &	     50 	&	 40.07 	&	 1.941 	&	 \textbf{0.019} 	&	1.517	&	&	96.82	&	97.51	&	\textbf{98.20}\\
% &	     60 	&	 116.0 	&	 5.332 	&	 \textbf{0.041} 	&	4.872	&	&	96.80	&	97.91	&	\textbf{98.47}\\
% &	     70 	&	 283.5 	&	 12.68 	&	 \textbf{0.054} 	&	10.72	&	&	97.08	&	98.05	&	\textbf{98.61}\\
% &	     80 	&	 440.0 	&	 20.77 	&	 \textbf{0.079} 	&	17.49	&	&	97.65	&	98.34	&	\textbf{98.61}\\
% &	     90 	&	 833.9 	&	 42.17 	&	 \textbf{0.102} 	&	30.04	&	&	97.36	&	98.33	&	\textbf{98.61}\\
% &	     100 	&	 2478.4 	&	 82.54 	&	 \textbf{0.129} 	&	66.97	&	&	97.22	&	98.47	&	\textbf{98.75}\\

%     \bottomrule
% \end{tabular}
% }
% % \label{tb:KPcpu}
% \end{table*}

\begin{table*}[!ht]
% \centering
% \begin{sidewaystable}
\centering
\makebox[\linewidth]{
 \begin{tabular}{cccccccccc} \toprule 
    &&&&\multicolumn{6}{c}{$PR$ variants} \\ \cline{5-10} 
    n&  KS* & FPBH & $RD$ & $PRrand$ & $PR$sim& $PR$dif & $PI$ &$PI$sim &$PI$dif\\ \midrule
    	10	&	6.58	&	\textbf{6.28	(95.4)}	&	5.91	(89.8)	&	6.02	(91.5)	&	5.96	(90.6)	&	5.95	(90.4)	&	6.03	(91.6)	&	5.94	(90.3)	&	5.93	(90.1)\\
    	20	&	6.97	&	6.70	(96.1)	&	6.61	(94.8)	&	6.75	(96.8)	&	6.72	(96.4)	&	6.71	(96.3)	&	\textbf{6.77	(97.1)}	&	6.71	(96.3)	&	6.70	(96.1)\\
    	30	&	7.21	&	6.92	(96.0)	&	6.96	(96.5)	&	7.05	(97.8)	&	7.04	(97.6)	&	7.04	(97.6)	&	\textbf{7.06	(97.9)}	&	7.04	(97.6)	&	7.04	(97.6)\\
    	40	&	7.14	&	6.86	(96.1)	&	6.95	(97.3)	&	7.01	(98.2)	&	7.01	(98.2)	&	7.01	(98.2)	&	\textbf{7.03	(98.5)}	&	7.01	(98.2)	&	7.01	(98.2)\\
    	50	&	7.23	&	7.00	(96.8)	&	7.05	(97.5)	&	7.09	(98.1)	&	7.09	(98.1)	&	7.09	(98.1)	&	\textbf{7.10	(98.2)}	&	7.09	(98.1)	&	7.09	(98.1)\\
    	60	&	7.18	&	6.95	(96.8)	&	7.03	(97.9)	&	7.06	(98.3)	&	7.06	(98.3)	&	7.06	(98.3)	&	\textbf{7.07	(98.5)}	&	7.06	(98.3)	&	7.06	(98.3)\\
    	70	&	7.19	&	6.98	(97.1)	&	7.05	(98.1)	&	7.08	(98.5)	&	7.08	(98.5)	&	7.08	(98.5)	&	\textbf{7.09	(98.6)}	&	7.08	(98.5)	&	7.08	(98.5)\\
    	80	&	7.22	&	7.05	(97.6)	&	7.10	(98.3)	&	\textbf{7.12	(98.6)}	&	7.11	(98.5)	&	\textbf{7.12	(98.6)}	&	\textbf{7.12	(98.6)}	&	7.11	(98.5)	&	\textbf{7.12	(98.6)}\\
    	90	&	7.20	&	7.01	(97.4)	&	7.08	(98.3)	&	\textbf{7.10	(98.6)}	&	\textbf{7.10	(98.6)}	&	\textbf{7.10	(98.6)}	&	\textbf{7.10	(98.7)}	&	\textbf{7.10	(98.6)}	&	\textbf{7.10	(98.6})\\
    	100	&	7.19	&	6.99	(97.2)	&	7.08	(98.5)	&	7.09	(98.6)	&	7.09	(98.6)	&	7.09	(98.6)	&	\textbf{7.10	(98.7)}	&	7.09	(98.6)	&	7.09	(98.6)\\
    \bottomrule
\end{tabular}
}
\caption{Comparing HV(\%) of algorithms on MOKP for $p$=3, * indicates optimal values, best heuristic values are in bold.} 
\label{tb:KPhv}
\end{table*}
% \end{sidewaystable}

In general, both the number of solutions and computation time increase as the size of instances becomes larger. The difference between FPBH and \textit{Bensolve} is clearly noticeable in Table \ref{tb:AP}, which shows that $Bensolve$ outperforms FPBH in all the MOAP instances. Furthermore, the difference becomes greater as the size of the instances grows. For example, for the instances with more than 15 tasks (n$\geq$15), the number of solutions of $Bensolve$ is more than double that of FPBH. Further, solutions are found in significantly less time.
Hence, not only the quantity but also the quality of the solutions of $Bensolve$ are better than FPBH. It reaches more than 99\% of the maximum HV throughout all the problem sets. Furthermore, the HV value increases as the problem size gets larger. On the other hand, the highest HV value of FPBH is 97.49\% in the smallest problem class (n=5), and it decreases as the problem size increases. This suggests that MOAP is not a suitable benchmark problem for MOIP heuristics.

In the case of the KP, FPBH and the $PR$ variants are competitive in terms of the number of solutions for instances with fewer than or equal to 70 items (n$\leq$ 70). For $n \geq$ 80, FPBH generates more solutions than $PR$ variants. 
$PR$ variants take less computation time than FPBH in all instances. Notably, $PRrand$ takes less than a quarter of the time than FPBH does. Skipping \textit{SelectionRules} and the best move analysis, but instead relying on randomness helps to reduce CPU time.
Except for $PI$, embedding the $ImprovedND$ operation into the $PR$ heuristics does not bring any noticeable difference in the number of solutions or run time. $PI$ finds the most solutions among all $PR$ variants while its CPU time increases fivefold. We also observe that it finds better solutions while spending longer running time in Table \ref{tb:KPhv}. 
$RD$ always produces the fewest solutions among the heuristic methods. However, it takes considerably less CPU time. For instance, it takes less than 1 second regardless of the problem size as it is seen in Table \ref{tb:KPcpu}. 
Although FPBH finds more solutions for larger instances, %these solutions are not necessarily good. In fact, 
every $PR$ variant achieves a higher HV than FPBH, except for the smallest problem class with n=10. In particular, $PI$ generates the highest quality solutions throughout the experiments. $RD$ also outperforms FPBH for the instances with more than 30 items.

We observe that the $RD$ and $PR$ variants do not show better performance than FPBH on the smallest instances. The reason for this may be the structure of the test instances.
When a problem size is small, a small fractional value has a relatively big impact on each objective. 
When a fractional value is rounded down, the solution quality deteriorates comparably more on smaller instances. In addition, the number of initially provided LB set solutions is limited in smaller instances. For these reasons, $RD$ has a large HV gap.
The quality of $RD$ also influences that of $PR$ variants. If very small IR sets are provided, 
the choice of the $S_i$-$S_g$ pair is restricted. This causes a limited number of new paths to be generated. 
% This means that after the \textit{rounding down} operation, we deduct the contribution of a fractional value which leads to a significant difference in a small instance.
% Also $RD$ has a large HV gap. The $PR$ variants are influenced by the number of solutions of $RD$.
% , which is not sufficient to find new solutions.
% the ability to create the pair of an initial and a guiding solution is limited. Thus, the algorithm cannot generate new paths.\\  

\section{\uppercase{Conclusion}} \label{sec::conclusion}
In this study, we propose a LP relaxation-based matheuristic for three-objective binary IP. The proposed algorithm relies on a high-performing vector LP solver, $Bensolve$, which provides bound sets, \textit{rounding down} and $PR$. In the computational study, we show that simple \textit{rounding down} can already find high-quality solutions in most instances. 
After embedding $PR$ with the first stage, the proposed heuristic generates more solutions, which show higher quality than that of FPBH in most problem classes. The number of solutions and CPU times of the $PR$ variants are similar to each other. Notably, the $PRrand$ algorithm takes less than a quarter of the computation time FPBH does. The biggest advantage of the proposed algorithm is that it can find high-quality approximations fast, which also shows its effectiveness.
For future work, we plan to extend the algorithm to deal with general MOIPs. 
% Further, it could be tailored to real-world applications such as supply chain network design. As the size of real-world problems be much larger, approaches to further reduce the computation time will be investigated. 
% In addition, the proposed method can be embedded into an interactive algorithm in order to facilitate decision making.

%Different ways of refining the best move analysis can be studied to find solutions from the unexplored search region.

% \section{\uppercase{Copyright Form}}

% \noindent For the mutual benefit and protection of Authors and
% Publishers, it is necessary that Authors provide formal written
% Consent to Publish and Transfer of Copyright before publication of
% the Book. The signed Consent ensures that the publisher has the
% Author's authorization to publish the Contribution.

% The copyright form is located on the authors' reserved area.

% The form should be completed and signed by one author on
% behalf of all the other authors.

% ADD AFTER REVIEWING PROCESS
\section*{\uppercase{Acknowledgements}}
This research was funded in whole, or in part, by the Austrian Science Fund (FWF) [P 31366-NBL]. For the purpose of open access, the author has applied a CC BY public copyright licence to any Author Accepted Manuscript version arising from this submission.

\bibliographystyle{apalike}
{\small
\bibliography{references}}

\section*{\uppercase{Appendix}}
Table \ref{tb:solvers} shows the comparison of three state-of-the-art algorithms that can deal with multi-objective linear programming.
\textit{Bensolve} (Ben) by \citet{lohne2017vector} and \textit{Inner solver} (Inner) by \citet{csirmaz2020inner} are
implemented in C and
publicly available at \url{http://www.bensolve.org/} and \url{https://github.com/lcsirmaz/inner}, respectively.
The algorithm suggested by \citet{ozpeynirci2010exact} (ÖK) is implemented in Julia.
GLPK is used as LP solver. The time limit of the experiment is 3600 seconds.
All the experiments are carried out on a Quad-core X5570 Xeon CPUs @2.93GHz with 48GB RAM.
The figures are the average results of 10 test instances over 10 runs.
As benchmark instances, we used the multi-objective assignment problem (MOAP), the multi-objective knapsack problem (MOKP), and multi-objective general integer linear programming problems (MOILP) which are all generated by \citet{kirlik2014new} and available at \url{http://home.ku.edu.tr/~moolibrary/}.
Each problem class is divided into subclasses. The subclasses are categorised by the number of items. For the MOAP, it is categorised by 5/10/15/30/50, whereas for the MOKP and MOILP, the subclasses are 10/30/50/70/100. Each subclass has 10 instances; thus, there are 50 instances per class in total.
We report the CPU time (sec) and the number of LPs. % solved during all the iterations.\\
n/a. indicates that the algorithm %could not complete finding a solution to all the instances 
did not terminate within the time limit. $^{*number}$ indicates the number of instances solved out of 10 instances.

Throughout the experiment, Ben and Inner are highly competitive. In terms of the number of solved LPs, the \textit{Inner solver} comprehensively outperforms the other two methods.
Overall, the \textit{Inner solver} performs the best. However, it does not provide all the bound set information we need for our heuristic. % information. 
Thus, we use $Bensolve$ to obtain the initial bound sets.

% \begin{sidewaystable}[b!]
\newpage
\begin{table}[!htbp]
 \centering
% \captionsetup{justification=centering}
% \rotatebox{90}{\begin{varwidth}{\textheight}\centering
 \begin{tabular}{crrrrrrrr}
 \toprule 
  &&\multicolumn{3}{c}{CPUtime(sec)} & & \multicolumn{3}{c}{\#solved LP} \\ 
  \cline{3-5} \cline{7-9} 
    Problem& \#item &Ben & Inner & ÖK& &Ben & Inner & ÖK \\ \midrule
    \multicolumn{1}{c}{\multirow{5}{*}{MOAP}} &5  & 0.004  & \textbf{0.003}& 2.46 && 30.2 & \textbf{23.0} & 42.0 \\ 
    &10  & 0.03  & \textbf{0.02} &30.20 && 117.2 & \textbf{110.1} &1916.8\\  
    &15  & 0.10 & \textbf{0.07} &844.04 &&237.6 & \textbf{230.5} &12845.7\\ 
    &30  & \textbf{1.57}&1.59& n/a. && 1015 & \textbf{1008} &n/a.\\
    &50  & \textbf{13.51} & 19.09 & n/a. && 2705 & \textbf{2698} & n/a.\\
    \\
    \multicolumn{1}{c}{\multirow{5}{*}{MOKP}} & 10  & 0.003 &\textbf{0.002} &4.27 & &46.0 &\textbf{39.0} &300.7\\
    &30  & 0.02&  \textbf{0.01} &320.98 & & 223.1&\textbf{216.0} &7749.0\\ 
    &50  & 0.03 &\textbf{0.02} &320.98 &&464.6 & \textbf{454.7}&n/a.\\
    &70  & \textbf{0.07} &\textbf{0.07} &n/a. && 978.6 & \textbf{920.0}&n/a.\\
    &100  & 0.16 &\textbf{0.14} &n/a. && 1962.0 & \textbf{1397.0}&n/a.\\
    \\
    \multicolumn{1}{c}{\multirow{5}{*}{MOILP}} & 10  & 0.003 & \textbf{0.001} &$3.71^{*7}$  &&  28.1& \textbf{18.3} &$12.97^{*7}$\\
    & 30  & \textbf{0.01} & 0.05 &$449.9^{*8}$  && 79.5& \textbf{70.1} &$3481.8^{*8}$\\
    & 50  &0.02& \textbf{0.01}  &$1295.4^{*2}$  &&134.71& \textbf{124.3} &$3764.5^{*2}$\\
    & 70  &0.04& \textbf{0.02} & $32.1^{*2}$  && 178.2& \textbf{167.0} &$857.0^{*2}$\\
    & 100  & 0.07&\textbf{0.03} &$3082.2^{*1}$  && 237.4& \textbf{223.53} &$12835.0^{*1}$\\
    \bottomrule
\end{tabular}
% \end{varwidth}
% }
\caption{Comparing algorithms on MOLP instances, $p$=3} \label{tb:solvers}
\end{table}
% \end{sidewaystable}

\end{document}